\documentclass[12pt,reqno]{amsart}
\usepackage{amsthm}

\usepackage{amssymb}
\usepackage[all]{xy}
\usepackage{hyperref}

\newtheorem{theorem}{Theorem}[section]

\theoremstyle{definition}

\numberwithin{equation}{section}

\begin{document}

\baselineskip=15.5pt

\title[Lie's classification of nonsolvable subalgebras of vector fields]{On Lie's
classification of nonsolvable subalgebras of vector fields on the plane}

\author[H. Azad]{Hassan Azad}

\address{Abdus Salam School of Mathematical Sciences, GCU, Lahore 54600, Pakistan}

\email{hassan.azad@sms.edu.pk}

\author[I. Biswas]{Indranil Biswas}

\address{Department of Mathematics, Shiv Nadar University, NH91, Tehsil Dadri,
Greater Noida, Uttar Pradesh 201314, India}

\email{indranil.biswas@snu.edu.in, indranil29@gmail.com}

\author[M. A. Fazil]{M. Ahsan Fazil}

\address{Abdus Salam School of Mathematical Sciences, GCU, Lahore 54600, Pakistan}

\email{ahsan.fzl@sms.edu.pk}

\author[F. M. Mahomed]{Fazal M. Mahomed}

\address{DSI-NRF Centre of Excellence in Mathematical and Statistical Sciences, School
of Computer Science and Applied Mathematics, University of the Witwatersrand,
Johannesburg, Wits 2050, South Africa}

\email{Fazal.Mahomed@wits.edu.za}

\subjclass{17B10, 17B66, 32M25}

\keywords{Root systems, vector fields, foliation, Levi decomposition, highest weight vectors}

\date{}

\begin{abstract}
A brief proof of Lie's classification of finite dimensional subalgebras of vector
fields over the complex plane that have a proper Levi decomposition is given. The proof uses
representation theory of $sl(2,{\mathbb C})$. This, combined with \cite{ABF2} and
\cite{ABF3} completes the classification of finite dimensional subalgebras of vector
fields over the complex plane.
\end{abstract}

\maketitle

\section{Introduction}\label{se1}

Lie's results in his collected papers \cite[Vol.~3]{LiEn1} and their proofs continue to remain difficult to 
access. The main reason seems to be the sheer volume of his work, the language and above all the ideas that 
are routinely applied by contemporary Lie theorists --- like those of highest weights and root systems --- 
were invented after his time. If one applies these ideas, the proofs simplify considerably --- as shown in 
this paper.

Lie always worked with complex analytic vector fields and in this paper all Lie algebras considered are of the 
same type.

The classification of Lie is local: he considered two Lie algebras of vector fields $L_1$ and $L_2$ defined on 
connected open subsets $U_1\, \subset\, {\mathbb C}^N$ and $U_2\, \subset\, {\mathbb C}^N$ respectively to be 
equivalent if by a local analytic change of variables, one can be transformed onto the other; see \cite{ABF3} 
for a more precise formulation.

Lie's classifications were based on the nature of generic orbits, whether the algebra 
was transitive or intransitive, whether it was primitive or imprimitive --- meaning that it had an invariant 
foliation (see Section \ref{se2} for a precise definition). The proofs of Lie are given in 
\cite[p.~1--109]{LiEn1} and in the English translation \cite[p.~40--131]{Me}.

The classifications are not 
given in terms of the algebras being solvable, semisimple or Levi--decomposable, but in geometric terms. The 
book \cite{AH} is also a translation of \cite{LiEn1} with commentaries by the authors.

The paper of 
Gonz\'alez-L\'opez, Kamran and Olver, \cite{GKO}, is concerned mostly with primitive Lie algebras on the real 
plane and they do not give details about the imprimitive Lie algebras, because they prove that the 
classification of imprimitive Lie algebras in the complex and real plane are the same and refer to 
\cite{LiEn1} for proofs.

In the paper \cite{ABF2}, the authors reprove the results of Lie for solvable Lie 
algebras of vector fields over the plane using representation theory and partial differential equations. The 
paper \cite{ABF3} classifies semisimple Lie algebras of vector fields over ${\mathbb C}^N$ of maximal rank using 
representation theory and the theory of root systems, as well as the results and methodology of \cite{ABF1}, 
which is a basic paper regarding the algebraic viewpoint on the classifications of finite dimensional semisimple 
algebras of vector fields on ${\mathbb C}^N$.

In this paper, using similar ideas, we give details for the 
classification of finite dimensional Lie algebras of vector fields on the plane that have a proper Levi 
decomposition.

We use the same equivalence relation on Lie algebras of vector fields as mentioned above.
Two Lie algebras of vector fields $L_1$ and $L_2$ defined on connected
open subsets $U_1\, \subset\, {\mathbb C}^2$ and $U_2\, \subset\, {\mathbb C}^2$ respectively
will be called equivalent if there are connected open subsets $U\, \subset\,U_1$ and $V\, \subset\, U_2$
together with an analytic diffeomorphism
$$
\Phi\ :\ U\ \longrightarrow\ V
$$
such that $\Phi_* L_1\ =\ L_2$ on $V$. The above mentioned classification is up to this
equivalence.

The main result of this paper is as follows:

\begin{theorem}\label{t1}
If $\mathfrak g$ is a finite dimensional Lie algebra of vector fields on the plane that has a proper Levi 
decomposition $${\mathfrak g}\ =\ \mathfrak{s} \ltimes \mathfrak{r},$$ then $\mathfrak{s}$ must be isomorphic 
to $sl(2,{\mathbb C})$ and either $\mathfrak{r}$ is abelian or it has an abelian ideal of codimension $1$.

The Lie algebra $\mathfrak{r}$ produces a $1$--dimensional invariant foliation given by $\langle 
\partial_y\,\rangle$ or $\langle\, \partial_x\, +\,y\,\partial_y\,\rangle$ in suitable coordinates.
\end{theorem}

A more precise formulation of Theorem \ref{t1} is given in Theorem \ref{thm1}. This result can also be deduced 
by inspection of the tables in \cite[p.~3]{GKO}, \cite[p.~369--372]{AH} and \cite[p.~129--131]{Me}. We have 
derived the result by using the representation theory of $sl(2,{\mathbb C})$.

\section{Sketch of the argument}\label{se2}

Before giving a sketch of the argument, let us put down some definitions and notation that will be used in the 
rest of the paper.

If $U$ is an open subset of ${\mathbb C}^N$, by $V(U)$ we mean the space of all 
analytic vector fields defined on $U$. Following Lie, we often write this as $V({\mathbb C}^N)$.

An irreducible representation is determined by its highest weight vector. The
highest weight vector is unique up to multiplication by a nonzero scalar. If $V$ is a highest weight vector of an 
irreducible representation, so is $cV$ for all nonzero constant scalars $c$. We will refer to $V$ as a
highest weight vector of the irreducible representation.

The following formula will be used repeatedly 
without further comments: $$ [\exp(\chi)u,\, \exp(\psi)v]\ =\ \exp(\chi+\psi)(u(\psi)v-v(\chi)u +[u, \,v]); $$ 
here $u$ and $v$ are vector fields and $\chi$, $\psi$ are differentiable functions.

A foliation on a smooth manifold $M$ is $C^\infty$ subbundle ${\mathcal F}\ \subset\ TM$ such that the
sheaf of sections of $\mathcal F$ are closed under the operation of Lie bracket of vector fields. The
dimension of a foliation defined by $\mathcal F$
is the rank of the vector bundle ${\mathcal F}$. If $\dim M\,=\, d$ and
$\mathcal F$ is a foliation on $M$ of dimension $d'$, then for every point $x\, \in\, M$, there is an
open neighbourhood $x\, \in\, U_x\, \subset\, M$ together with a $C^\infty$ submersion $\varphi\, :\,
U_x \, \longrightarrow\, {\mathbb R}^{d-d'}$ such that the subbundle ${\mathcal F}\big\vert_{U_x}
\, \subset\, TU_x$ coincides with the kernel of the differential
$$
d\varphi\ :\ TU_x\ \longrightarrow\ \varphi^*T{\mathbb R}^{d-d'}.
$$ 
An identical result holds in the holomorphic category. See \cite{Su} for foliation associated to algebras
of vector fields.

If a foliation is 
generated by a subalgebra $\mathfrak h$ of vector fields on a manifold $M$, then it is $\mathfrak 
g$--invariant if and only if $$ [{\mathfrak g},\, {\mathfrak h}]\ \subset\ {\mathfrak h}. $$ If $L$ is a
Lie algebra of 
vector fields defined on an open subset $U$ of ${\mathbb C}^N$, then the rank of $L$ is, by definition $$ {\rm 
Max}_{p\in U} \dim \{X(p)\,\, \big\vert\,\,\, X\, \in\, L\}. $$ We use the fact that the only semisimple 
algebra of vector fields on ${\mathbb C}^2$ that has a solvable extension as vector fields on ${\mathbb C}^2$ 
is $sl(2,{\mathbb C})$ \cite[Corollary 4.1]{ABF3}. There are two types of realizations of $sl(2,{\mathbb C})$ 
as vector fields on ${\mathbb C}^2$:
\begin{enumerate}
\item[(t(1))] Those with generators $X\,=\, 
\exp(x)\partial_x$,\, $Y\,=\, - \frac{\exp(-x)\partial_x}{2}$.

\item[(t(2))] Those with generators $$X\,=\, 
\exp(x)\partial_y,\ \ \, Y\,=\, \exp(-x)(y\partial_x +(\frac{y^2}{2}+\varepsilon)\partial_y),$$ where 
$\varepsilon\,=\, 0$ or $1$.
\end{enumerate}
We want to determine extensions ${\mathfrak g}\,=\, sl(2,{\mathbb 
C})\ltimes {\mathfrak{r}}$, where $\mathfrak{r}$ is a solvable ideal of $\mathfrak g$. The radical 
$\mathfrak{r}$ decomposes into irreducible representations of $sl(2,{\mathbb C})$ and each irreducible 
representation is determined by its highest weight vector \cite[p.~59]{Ki} (recall that the highest
weight vector is unique up to multiplication by a nonzero scalar).

For Lie algebras of type t(1), all 
irreducible representations are one dimensional having highest weight $f(y)\partial_y$. Therefore, 
$\mathfrak{r}$ is a solvable algebra of vector fields on a line.

For type t(2), we note that every solvable 
extension contains an abelian extension. If $V$ is a highest weight vector in this abelian extension, then $$ 
[V,\, [Y,\, V]]\ =\ 0. $$ Using this equation, we find that the Lie algebras of type t(2) are extendable only if 
$\varepsilon \,=\, 0$. We then determine the highest weight vectors in any solvable extension whose weight --- 
meaning the eigenvalue of $H\,=\, [X,\, Y]$ --- is as large as possible. There are only two types of such 
weights, namely $$ V_1\ =\ \exp(x/2)\partial_y \ \ \, \text{ and }\ \ \, V_2\ =\ 
\exp(dx)(\partial_x+y\partial_y).$$ Working with the canonical coordinates for $Y$ we find that the radical 
$\mathfrak{r}$ is either abelian or it has an abelian ideal of codimension one.

In the next section, we 
prove the main result, following this outline. \section{Statement and proof of the main result}\label{se3} In 
the following statements, if $V_k$, $k\,=\,1,\,2,\, \cdots,$ are vectors in a representation space of a Lie 
algebra $\mathfrak g$, by $$ \langle V_1,\, V_2,\, \cdots\rangle $$ we mean the $\mathfrak g$--invariant 
subspace of the representation space generated by $V_1,\, V_2,\, \cdots$.

\begin{theorem}\label{thm1}\mbox{}
\begin{enumerate}
\item[{\textbf{(a):}}]\, If a semisimple algebra $\mathfrak{s}$ has a non-trivial solvable extension in V$ 
(\mathbb{C}^2)$ then $\mathfrak{s}$ must be isomorphic to $sl(2,{\mathbb C}).$

\item[{\textbf{(b):}}]\, 
The only realization of $sl(2,{\mathbb C})$ in V$(\mathbb{C}^2)$ are given --- up to the equivalence relation
on the Lie algebras of vector fields --- by $$(I):\,\,
X\,=\,\exp(x)\,\partial_x, \ \ \ Y\,=\, -\frac{\exp(-x)\,\partial_x}{2}, \ \ \ H\,=\,[X,\, Y].$$
$$(II):\,\, X\ =\ \exp(x)\,\partial_y,\, \ \ \, Y\ =\ \exp(-x)(y\partial_x +(\frac{y^2}{2}+\varepsilon)\partial_y),
\ (\varepsilon\ =\ 0,\,1),$$
$$
H\ =\ [X,\, Y].$$

\item[{\textbf{(c):}}]\, If $sl(2,{\mathbb C})$ is of type $\textbf{(b)}(I)$ then its extensions are 
\begin{itemize}
\item $sl(2,{\mathbb C})\, \ltimes\,\langle\,\partial_y\,\rangle$, and

\item $sl(2,{\mathbb C})\, \ltimes\,\langle\,\partial_y,\,y\,\partial_y\rangle.$
\end{itemize}
An invariant foliation is 
given by $\langle\,\partial_y\,\rangle$.

\item[{\textbf{(d):}}]\, If $sl(2,{\mathbb 
C})$ is of type $\textbf{(b)}(II)$ then it has a solvable extension only if $\varepsilon\,=\, 0$. In this 
case, a highest weight vector of maximal weight in the radical $\mathfrak{r}$ is $\exp(x/2)\partial_y$ or 
$\exp(dx)(\partial_x+y\partial_y)$, where $2d$ is a nonnegative integer (recall that the highest
weight vector is unique up to multiplication by a nonzero scalar).

The radicals are: 
\begin{itemize}
\item $\langle \exp(x/2)\partial_y\rangle$,

\item $\langle\exp(x/2)\partial_y\rangle\rtimes 
\langle\partial_x+y\partial_y\rangle$,

\item $\langle\exp(dx)(\partial_x+y\partial_y)\rangle$,

\item $\exp(dx)(\partial_x+y\partial_y)\rtimes\langle\partial_x+y\partial_y\rangle$.
\end{itemize}
An invariant foliation is given by $\langle\,\partial_x\,+\, y\,\partial_y\rangle.$

\end{enumerate}
\end{theorem}

\begin{proof}
Recall that a highest weight vector of weight $d$ in a representation of $sl(2,{\mathbb C})$ is a vector $V$ 
such that $X\cdot V\,=\, 0$ and $H\cdot V\,=\, dV$; here $X$, $Y$, $H$ are generators of $sl(2,{\mathbb C})$ 
such that $[X,\, Y]\,=\, H$, $[H,\, X]\,=\, X$, $[H,\, Y]\,=\, -Y$. By representation theory of $sl(2,{\mathbb 
C})$ \cite[p.~59, 4.8]{Ki}, $2d$ is a non-negative integer.

Now a highest weight vector is a common eigenvector for a Cartan subalgebra. Here the Cartan subalgebra is 
 spanned by $\partial_x$. As the solutions of the equation $\frac{\partial}{\partial x} f(x,\,y)\ =\
\lambda f(x,\,y)$ are of the form $\exp (\lambda x) g(y)$ and a highest weight vector $V$ is an
eigenvector for $\partial_x$ with $[X,\,V]\ =\ 0$, we see that a highest weight vector for the
given realization of $sl(2,{\mathbb C})$ is
is of the form
$$
V\ =\ \exp(dx)(f(y)\,\partial_x +g(y)\,\partial_y).
$$
Parts (a) and (b) are in \cite{ABF3}. The main 
point is that for a semisimple subalgebra of vector fields over $\mathbb{C}^N$ the rank of
its Cartan subalgebra is at most $N$, and a semisimple algebra of vector fields of rank $N$ on
$\mathbb{C}^N$ must be a product of type $A$ algebras.

For $sl(2,{\mathbb C})$ of type t(1), where $X\,=\, \exp(x)\,\partial_x$, $Y\,=\, -
\frac{\exp(-x)\,\partial_x}{2}$, the condition $[X,\, V]\,=\, 0$ implies that
$$
(d-1)f\ =\ 0. \ \ \text{ and }\ \ dg\ =\ 0.
$$
Therefore, either 
$V\,=\, g(y)\,\partial_{y}$ or $V\,=\,\exp{(x)} (f(y)\partial_{x})$. If $V\,=\,\exp{(x)}(f(y)\partial_{x})$, 
then $f$ cannot be a constant. Consequently,
$$
[Y,\, V]\,=\,f(y)\,\partial_x \ \ \text{ and } \ \ 
(\text{ad}(f(y)\partial_{x}))^{n}(V) \ =\ \exp (x)(f(y)^{n+1}\,\partial{x}).
$$ This means that the 
representation space is infinite dimensional. Hence all highest weight vectors are of the form 
$g(y)\,\partial_{y}$. Therefore, $\mathfrak{r}$ is a solvable Lie algebra of vector fields on the line. 
It can be shown that $$\mathfrak{r}\ = \ \langle{\partial_{y}\rangle} \ \ \text{ or } \ \ \mathfrak{r}\ =\ 
\langle{\partial_{y},\,y\,\partial_{y}\rangle}.$$
This can be proved as follows: If $A$ is a nilpotent algebra of vector fields on the line and $W$ is a nonzero 
vector in the center of $A$, then if $y$ is a coordinate on the line in which $W$ is locally $\partial_y$, 
then $A$ consists of all multiples of $\partial_y$. Therefore, if $\mathcal G$ is a nonabelian solvable Lie 
algebra of vector fields on the line, its commutator must be abelian. If $[{\mathcal G},\, {\mathcal G}]$ in 
local coordinates is $\langle \partial_y\rangle$, its normalizer is the space of all vector fields of the form 
$(\alpha y+\beta)\partial_y$, where $\alpha$ and $\beta$ are constants.

Assume now that
$$sl(2,{\mathbb C})\ =\ \langle{\exp(x)\,\partial_{y},\,\exp(-x)\,(y\,\partial_{x}\,+\,
(\frac{y^2}{2}\,+\,\varepsilon)\,\partial_{y}})\,\rangle.$$ 
Let $$X\ =\ \exp{(x)}\,\partial_{y},\ \ Y\ =\ \exp(-x)\,(y\,\partial_{x}
\,+\,(\frac{y^2}{2}\,+\,\varepsilon)\,\partial_{y}),\, \ (\varepsilon\,=\, 0,\, 1).
$$
Thus we have $H\,=\, [X,\, Y]\,=\, 
\partial_x$. Since $V\ =\ \exp(dx)(\,f(y)\,\partial_{x}\,+\,g(y)\,\partial_{y})$, it follows
from the statement $[X, \, V]\,=\,0$ that $$V\ =\ \exp(dx)(\kappa\,\partial_{x}\,+
(\kappa\, y+\ell)\,\partial_{y}).$$ Following the sketch given in 
Section 1, we divide the rest of the proof into several steps.

\textbf{Step 1.}\, A type \textbf{(b)} $sl(2,{\mathbb C})$ is extendable only if $\varepsilon\,=\,0.$

\begin{proof}[{Proof of Step 1}]
Assume that we have an extension $sl(2,{\mathbb C})\ltimes {\mathfrak{r}}$, where $\mathfrak{r}$ is solvable. 
Thus we have $\mathfrak{r'} \, :=\, [\mathfrak{r},\, \mathfrak{r}]\, =\, 0$ or else $sl(2,{\mathbb 
C})\ltimes\, Z(\mathfrak{r'})$ gives an abelian extension, where $Z(\mathfrak{r'})$ is the center of 
$\mathfrak{r'}$.

This means that if we take a highest weight vector $V$ in an abelian extension $\mathfrak{r}$, then $[V, \, 
[Y, \, V]]\ =\ 0$.

Assume now that $V$ is a highest weight vector in this abelian extension of maximal weight $d$. We have $[V, 
\, [Y, \, V]]\ =\ 0$. If $\text{rank}\, \langle X,\, V\rangle\ =\ 1,$ we can take $V\ =\ 
\exp(dx)(\partial_{y})$. Then $$[V, \, [Y, \, V]]\ =\ \exp((2d-1)x)\,((2d-1)\,\partial_{y}).$$ This is $0$ if 
and only if $d\,=\,\frac{1}{2}$. In this case, the subspace generated by $V$ is $2-$dimensional, and therefore, 
$(\text{ad}(Y))^2(V)\, =\, 0$. Now $$(\text{ad}(Y))^2(V)\ =\ 
\exp{(-x(1+\frac{1}{2}}))\,(\varepsilon(1+ \frac{1}{2})\partial_{y}).$$ Consequently, $\varepsilon$ must be 
zero if $\text{rank}\, \langle X,\, V\rangle\ =\ 1.$

Assume now that $$\text{rank}\, \langle X,\, V\rangle\ =\ 2.$$ Thus $V$ is of the form 
$$\exp(dx)\,(\partial_x+(y+\ell)\,\partial_y).$$ The $\partial_x$ component of $[V, \, [Y, \, 
V]\,]$ is, in general, $\exp (x(2d-1))\,(\ell\,(d+1))$. Therefore, $\ell\ =\ 0$ and the 
$\partial_y$ component of $[V, \, [Y, \, V]]$ is $\exp (x(2d-1))\,(-4\epsilon\,+2d\, 
\varepsilon))$. If $\varepsilon\,\neq\, 0$ then $d$ must be $2$.

Now we check whether $V\ =\ \exp{(2x)}\,(\partial_x+\, y\partial_y)$ gives an abelian extension.

Now $\text{ad}(Y)^2(V)$ is of weight $0$ and $$\mathcal{U}\ =\ \text{ad}(Y)^2(V)\ =\
3\,y^{2}\,\partial_x\,+(3\,y^3\,+6\,y)\partial_y.$$ Note that $$[\mathcal{U},\, \text{ad}(Y)^3(V)]\ =\ [\mathcal{U}, \, \exp(-x)(3y^3\,\partial_x\,+3(y^4\, +3y^2\,+2))\,\partial_y]\ \neq\, 0$$
as seen by considering its $\partial_x$ component --- which is $\exp(-x)y$ (up to a nonzero
scalar multiple). Thus, the extension is not abelian, hence $\varepsilon$
must be $0$ for an $sl(2,{\mathbb C})$ to have a solvable extension. This establishes Step 1.
\end{proof}

\textbf{Step 2.}\, If $\varepsilon\ =\ 0$ and we have a solvable extension
$sl(2,{\mathbb C})\ltimes {\mathfrak{r}}$ and $V$ is a highest weight vector of maximal weight $d$
with $\text{rank}\, \langle X,\, V\rangle\ =\ 2$, then $V\ =\ \exp(dx)(\partial_x\,+y\,\partial_y)$
up to a nonzero scalar multiple (recall that the highest weight vector is unique up to
multiplication by a nonzero scalar).

\begin{proof}[Proof of Step 2]
As $\text{rank}\, \langle X,\, V\rangle\ =\ 2$, $V$ is of the form 
$$ V\ =\ \exp(dx)(\partial_x\,+(y\,+\,\ell)\,\partial_y).$$
Now $$[V, \, [Y, \, V]\,]\ =\ \exp(x(2d-1))\,(\ell(d+1)\,\partial_x\,+(\ell^2(2d-1)\,+\,
\ell y\,(d+1)\,\partial_y)). $$
For $d\,>\,1,$ we have $$[V, \, [Y, \, V]\,]\ =\ 0,$$ because if it is not zero, it would be weight $2d-1$ and
$2d-1\,>\,d$ if $d\,>\,1.$ Thus, we have $\ell\ =\ 0$ for $d\,>\,1.$

For $d\, =\, 1$ 
$$[V, \, [Y, \, V]\,]\ =\exp(x)\,(\,2\,\ell\,\partial_x\,+(\ell^2\,+\, 2\,y\,\ell)\,\partial_y\,) .$$

So, if $\ell\,\neq\,0,$ we see that 
$$\exp\,(x)(\,\partial_x\,+\,(\frac{\ell}{2}\,+\,y)\,\partial_y\,)\ \in\ \mathfrak{r}.$$ As 
$\exp\,(x)\,(\,\partial_x\,+\,(y\,+\ell\,)\,\partial_y)\,\in\,\mathfrak{r},$ it follows that 
$\exp(x)\,(\,\frac{\ell}{2}\,\partial_y\,)$ is in $\mathfrak{r}$ and therefore
$X\,=\,\exp(x)\,\partial_y$ is 
in $\mathfrak{r},$ thus leading to a contradiction.
Therefore, for $d\,=\,1,$ it follows that $\ell$ must be zero.

For $d\,=\,\frac{1}{2}$ and $\ell\,\not=\, 0$,
$$[V, \, [Y, \, V]\,]\ =\,\partial_x\,+\,y\,\partial_y\quad (\text{up to a non-zero multiple}).$$
The commutator of $\partial_x\,+\,y\,\partial_y$ with $V$ is
$\frac{1}{2}\exp(\frac{x}{2})(\,\partial_x\,+\,(y-\,\ell)\partial_y\,).$ Therefore, if $\ell\,\neq 0$ $\exp(\frac{x}{2})\,\partial_y\,\in\,\mathfrak{r}.$ Its commutator with $V$ is $\exp(x)\,(\,\frac{1}{2}\,\partial_y\,).$ Thus $X\,=\,\exp(x)\,\partial_y$ would be in $\mathfrak{r}.$

Finally, if $d\,=\,0,$ then $[\,Y\, ,V\,]$ must be zero. But $$[\,Y\, ,V\,]\,=\exp(-x)\,(\,2\,\ell\,\partial_x\,+\,(\,2\,y\,\ell\,)\,\partial_y).$$ Hence $\ell$ must be zero. Therefore, $\ell\ =\ 0$ in every case.
This proves Step 2.
\end{proof} 

\textbf{Step 3.}\, If V is as in Step $2$ and its weight $d\,>\,0,$ then $V$ is the only highest 
weight vector in $\mathfrak{r},$ except possibly for the vector $(\partial_x\, +\, y\,\partial_y)$
(recall that the highest weight vector is unique up to multiplication by a nonzero scalar).

Moreover, $\mathfrak{r}\ =\ \langle\,\exp(dx)\,(\partial_x\, +\, y\, \partial_y)\,\rangle$ or $\mathfrak{r}\ =\ \langle\,\exp(dx)\,(\partial_x\, +\, y\, \partial_y)\,\rangle\,\rtimes\, \langle\partial_x\, +\, y\, \partial_y\rangle.$

\begin{proof}[Proof of Step 3]
Let $W\,=\,\exp(\widetilde{d}x)\,(\partial_x\,+\,(y+\ell)\partial_y)$. Now $$[\,W\,,V\,]\,=\,\exp((\widetilde{d}\,+\,d\,)(x))\,(d\,(\,\partial_x\,+\,y\,
\partial_y\,)\,-\widetilde{d}\,(\,\partial_x\,+\,(\,y\,+\,\ell\,)\,\partial_y\,)\,+\,\ell\,\partial_y\,).$$
Now assume that $\widetilde{d}\,>\,0.$ Thus $[V,\,W],$ if it is not zero, is of weight $d+\widetilde{d}.$
Therefore, $[V,\,W]\,=\,0$ implies that $d\,=\,\widetilde{d}$ and $\ell\,(1-d)=\,0.$ As $d\,>\,1,$ we
must have $\ell\,=\,0.$ Thus we have $W\,=\,V.$ Finally, if $W\,=\,\exp(\widetilde{d}x)\,(\partial_y)$ and
$\widetilde{d}\,>\,0$ then $[V,\,W]\,=\,0$ 
implies that $\widetilde{d}\,=\,1,$ so $W$ would equal $X.$ Moreover, $\widetilde{d}\,=\,0$ would mean
that $[\partial_y,\,Y]\,=\,0.$ This establishes Step 3.
\end{proof}

\textbf{Step 4.}\, If $V$ is a highest weight vector of maximal weight in a solvable extension of type
\textbf{(b)} of $sl(2,{\mathbb C})$ with $\text{rank}\, \langle X,\, V\rangle\ =\ 1$
(recall that the highest weight vector is unique up to multiplication by a nonzero scalar), then $V\,=\,
\exp(\frac{x}{2})\,\partial_y$. The only other highest weight vector can be $(\,\partial_x\,+\,y\,\partial_y\,).$

\begin{proof}[Proof of Step 4]
We have $V\,=\,\exp(dx)\,\partial_y.$ As in Step $1$
$$[V,\,[Y,\,V]]\ =\ \exp((2d-1)x)\,((2d-1)\,\partial_y).$$
Now if $d\,>\,1,$ then $[V,\,[Y,\,V]]$ must be $0,$ as explained in Step $2$;
note that $d$ can't be $1$ or $0$ as in Step 3. Thus $d$ must be $\frac{1}{2}$ and
the irreducible subspace generated by $V$ is two dimensional and it is abelian. The only other highest
weight vector is $$W\ =\ (\partial_x\,+\,(y+\ell)\,\partial_y)$$
(recall that the highest weight vector is unique up to multiplication by a nonzero scalar).
Thus $[Y,\,W]$ must be zero and this
gives $\ell\,=\,0.$ This proves Step 4.
\end{proof}

\textbf{Step 5.}\, If $V$ is a highest weight vector of maximal weight in a solvable extension of 
$sl(2,{\mathbb C})$ and $\text{rank}\, \langle X,\, V\rangle\ $ is $2,$ then the irreducible subspace 
generated by $V$ is abelian and the only other possible highest weight vector in $\mathfrak{r}$ is 
$\partial_x\,+\,y\,\partial_y.$

\begin{proof}[Proof of Step 5]
We need only to prove that the subspace generated by $V$ is abelian, as the other assertions have already been proved in Step $3$. To show that the subspace generated by $V$ is abelian, we use coordinates $(\widetilde{x},\,
\widetilde{y})$ in which $Y\,=\,\partial_{\widetilde{x}}.$ We can take $$\widetilde{x}\,=\,
2\frac{\exp (x)}{y},\ \ \
\widetilde{y}\,=\,2\frac{\exp (x)}{y^2}.$$ In these coordinates, $$y\ =\
\frac{\widetilde{x}}{\widetilde{y}},\, \ \ \exp (x)\ =\ \frac{\widetilde{x^2}}{2\widetilde{y}},
\, \ \ \partial_x\ =\ \widetilde{x}\,\partial_{\widetilde{x}}\,+\,\widetilde{y}\,
\partial_{\widetilde{y}},\ \ \, y\partial_y\ =\ -\widetilde{x}\,\partial_{\widetilde{x}}\,-\,2\widetilde{y}\,\partial_{\widetilde{y}}.$$ 
Therefore, $$\partial_x\,+\,y\partial_y\,=\,-\widetilde{y}\partial_{\widetilde{y}}\ \ and\ \
V\,=\,-\frac{\widetilde{x}^{2d}}{2^{d}\widetilde{y}^{(d-1)}}\,\partial_{\widetilde{y}}.$$

Hence $$\langle\,\text{ad}(Y)^{n}(V)\,\,\big\vert\,\,\,n\in\,\mathbb{Z}\rangle\ =\ \langle\,
(\,\frac{{x}^{n}}{2^{d}\widetilde{y}^{(d-1)}}\,)\partial_{{y}}\,\,\big\vert\,\, 0 \,\leq\, n\,\leq\, 2d\rangle$$ and
this is clearly abelian. The subalgebras $\langle\,\partial_y\,\rangle$ and $\langle\,\partial_x\,+\,y\,
\partial_y\,\rangle$ are normalized by the algebras in \textbf{(b)}(I) and \textbf{(b)}(II), hence
they give an invariant foliation. This completes the proof of Step 5.
\end{proof}

Now the proof of the theorem is complete.
\end{proof}

\section*{Acknowledgements}

We thank the referee for helpful comments to improve the exposition.
IB is partially supported by a J. C. Bose Fellowship (JBR/2023/000003)


\end{document}